\documentclass[10pt,onecolumn,a4,reqno]{article}
\usepackage{amsmath, amsthm}
\usepackage{graphicx}
\usepackage{amssymb}
\usepackage{amsfonts}
\usepackage{delarray}
\usepackage{stackrel}
\usepackage{delarray}

\usepackage{amsmath, amsthm}

\title {\textbf{ A Proof of Sendov's  Conjecture }}
\author {G.M. Sofi }
\date{}
\begin{document}
\maketitle
\begin{center}
Department of Mathematics,\\
Central University of Kashmir,\\
Ganderbal-191201(India).\\
\end{center}
\begin{center}
{\bfseries Email:} { $gmsofi@cukashmir.ac.in $} 
\end{center}
\begin{abstract}
     Sendov's conjecture asserts that if all the zeros of a polynomial $p$ lie in the closed unit disk, then there must be a critical point of $p$ within unit distance of each zero. The conjecture has been proved to be true for many special cases (See [1]). Here is a  proof of the  conjecture for all polynomials. 
\end{abstract}
\maketitle{}
\section{Introduction}\label{sec.1}
\indent~~~ The Sendov's conjecture is an  unproven conjecture in complex analysis. This conjecture was included in the collection $Research ~Problems~in~Function~ Theory$ in 1967 by Professor Hayman. The conjecture is due to the Bulgarian Mathematician B. Sendov.  Till present the conjecture has been found to be true for many special cases(see[1]). What follows is a  proof of  the conjecture .\\

\noindent\textbf{Theorem 1}. Let $p(z) = z^n + a_{n-1}z^{n-1} + a_{n-2}z^{n-2} + \cdots + a_1 z + a_0$ be a non-constant polynomial with all its zeros $z_1, z_2, \cdots ,z_n$ lying inside the closed unit disk $|z|\leq 1$. Then each of the disks $|z-z_i|\leq 1,~i=1,2,\cdots , n$ must contain a zero of $p'$.\\\\
\textbf{Proof.} [we infact prove that if $ p $ has all its zeros $z_1, z_2, \cdots ,z_n$ lying inside any closed  disk $|z-a|\leq r$. Then each of the disks $|z-z_i|\leq r,~i=1,2,\cdots , n$ must contain a zero of $p'$ and then taking $ r=1$ and $a=0 $ we get the required conclusion.]\\~\\
 {\bf Case1:} $ \sum_{i=1}^{n}z_i=0 $ \\\\\
   ~~~~Let~$ \zeta _i,~i=1,2,\cdots ,n-1$ be the critical points of $p$. Assume to the contrary that the conclusion of the theorem does not hold. Then there exists a zero of $p$ say $z_1$ such that 
$
  |z_1-\zeta_i|>r~~\mbox {for all $1\leq i\leq n-1$}
$.
Therefore by Gauss-Lucas theorem all  $\zeta _is$ must  lie in the crescent-shaped region outside $|z-z_1|\leq r$ but inside $|z-a|\leq r$(as indicated in Figure 1). We can  draw a line through $ z_1 $ parallel to the line through the points   of intersection of the  circles $|z-a|\leq r$  and  $|z-z_1|\leq r$ (see Figure 1)  and it is clear that all the vectors $z_1 -\zeta _i,~i=1,2,\cdots ,n-1$~lie in the same half plane determined by  this line.Hence $ z_1 $ cannot be the centroid of these $z_1 -\zeta _is$.\\
\begin{figure}[h]
\centering
\includegraphics[scale=.6]{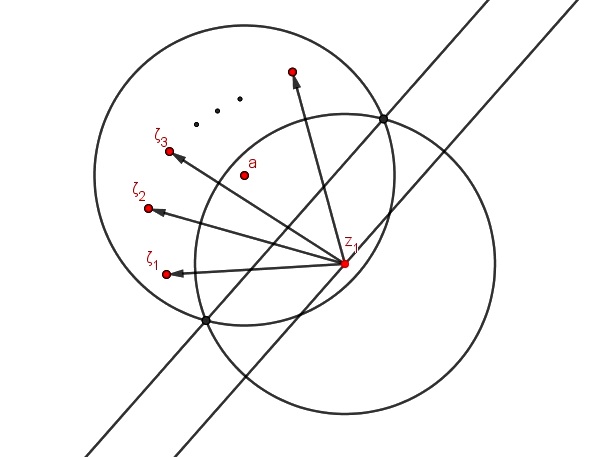}
\caption{$z_1 -\zeta _is$ are exactly the antiparallel vectors of the ones indicated.}
\end{figure}\\

 But  \\
\~~ \begin {equation*} z_1=z_1 -\frac{\sum_{i=1}^{n}z_i}{n}=z_1 -\frac{\sum_{i=1}^{n-1}\zeta _i}{n-1}= \frac{\sum_{i=1}^{n-1}(z_1-\zeta _i)}{n-1}\end {equation*}
\ a contradiction.
\\~\\
 {\bf Case2:} $ \sum_{i=1}^{n}z_i\neq 0 $ \\
 
 Now by case 1  the Sendov's conjecture holds for the polynomial $ p(z-\frac{a_{n-1}}{n }) $. Since a translation $ z\rightarrow z+\frac{a_{n-1}}{n } $ does not alter the relative distances between zeros and critical points of $ p(z-\frac{a_{n-1}}{n }) $ Sendov's conjecture must be true for $ p(z) $ also. 
 \qed

\end{document}